\documentclass[11pt]{article} \topmargin=-0.5cm
\usepackage{amssymb}

\newtheorem{theorem}{Theorem}
\newtheorem{lemma}{Lemma}
\newtheorem{corollary}{Corollary}
\newtheorem{condition}{Condition}
\newtheorem{remark}{Remark}
\def\e{\varepsilon}

\def\defi{\stackrel{{\scriptscriptstyle \Delta}}{=}}

\def\a{\alpha}
\def\d{\delta}

\def\w{\widehat}

\def\esssup{\mathop{\rm ess\, sup}}
\def\ulim{\mathop{\overline{\rm lim}}\,}

\def\R{{\bf R}}

\def\b{\beta}

\def\g{\gamma}

\def\ww{\widetilde}

\def\t{\theta}
\def\oo{\bar}

\def\p{\partial}

\def\A{{\cal A}}

\newcommand{\be}{\begin{equation}}
\newcommand{\ee}{\end{equation}}
\newcommand{\bd}{\begin{displaymath}}
\newcommand{\ed}{\end{displaymath}}
\newcommand{\ba}{\begin{array}{ll}}
\newcommand{\ea}{\end{array}}
\newcommand{\baa}{\begin{eqnarray}}
\newcommand{\eaa}{\end{eqnarray}}
\newcommand{\baaa}{\begin{eqnarray*}}
\newcommand{\eaaa}{\end{eqnarray*}}   \font\sm=cmr10

\def\ww{\tilde}

\def\CC{{\cal C}}

\date{ }
\title{Universal estimate of the  gradient for
parabolic equations
}
\author{
Nikolai Dokuchaev
\\ {\sm Department of Mathematics, Trent University, Ontario,
Canada}}
\begin{document}
\maketitle
\begin{abstract}
We suggest a modification of the  estimate for  weighted Sobolev
norms of solutions of parabolic equations such that the matrix of
the higher order coefficients  is included into the weight for the
gradient. More precisely, we found the upper limit estimate that can
be achieved by variations of the zero order coefficient.  As an
example of applications, an asymptotic estimate  was obtained for
the gradient at initial time.  The constant in the estimates is the
same for all possible choices of the dimension, domain, time
horizon, and the coefficients of the parabolic equation. As an
another example of application, existence and regularity results are
obtained for parabolic equations with time delay for the gradient.
\\ {\it AMS 2000 subject classification:}
35K10, 
35K15, 
35K20 
\\ {\it Key words and phrases:} parabolic
equations, regularity, solution gradient, asymptotic estimates,
parabolic delay equations.
\end{abstract}
\section{Introduction}
We study Dirichlet boundary value problems for linear parabolic
equations. The classical results give  estimates for Sobolev norms
of the solution via the $L_2$-norm of the nonhomogeneous  term
(see, e.g., the second energy inequality (4.56) in Ladyzhenskaia
(1985), Chapter III). We suggest a modification of this estimate
for weighted Sobolev norms  such that the matrix of the higher
order coefficients  is included into the weight for the gradient.
We found the limit upper estimate that can be achieved
 by variations of the zero order coefficient. More precisely, we
obtain estimates for  $e^{-Kt}u(x,t)$ via $L_2$-norm of
$e^{-Kt}h(x,t)$, where  $u(x,t)$ and $h(x,t)$ is the solution and
the nonhomogeneous term respectively, and where $K>0$ is being
variable. As an example of applications, an asymptotic upper
estimate was obtained for a weighted $L_2$-norm of the gradient at
initial time. The constants in these estimates are the same for
all possible choices of the dimension, domain, time horizon, and
the coefficients of the parabolic equations, i.e., these estimates
can be called universal estimates. As an another example of
applications, we establish solvability and regularity for special
parabolic equations such that the gradient is included with time
delay.
 \section{Definitions}\label{SecD}
\subsection*{Spaces and classes of functions.} 
 We denote by $|\cdot|$ the Euclidean norm in $\R^k$ and the Frobenius norm in $\R^{k\times
 m}$, and we denote by $\|\cdot\|_{ X}$ the norm in a linear normed space
 $X$. We denote by
 $(\cdot, \cdot )_{ X}$ the scalar product in  a Hilbert space $
X$. For a Banach space $X$, we denote by $C([a,b],X)$ the Banach
space of continuous functions $x:[a,b]\to X$.
\par
 We are given an open domain $D\subseteq\R^n$ such that either
$D=\R^n$ or $D$ is bounded  with $C^2$-smooth boundary $\p D$.
\par
Let $T>0$ be given, and let $Q\defi D\times (0,T)$.
\par
We denote by ${W_2^m}(D)$   the Sobolev  space of functions that
belong to $L_2(D)$ together with the distributional derivatives up
to the $m$th order, $m\ge 0$.
\par Let $H^0\defi L_2(D)$,
and let $H^1\defi \stackrel{\scriptscriptstyle 0}{W_2^1}(D)$ be
the closure in the ${W}_2^1(D)$-norm of the set of all smooth
functions $u:D\to\R$ such that  $u|_{\p D}\equiv 0$. Let
$H^2=W^2_2(D)\cap H^1$ be the space equipped with the norm of
$W_2^2(D)$. The spaces $H^k$ are Hilbert spaces, and $H^k$ is a
closed subspace of $W_2^k(D)$,  $k=1,2$.
\par
 Let $H^{-1}$ be the dual space to $H^{1}$, with the
norm $\| \,\cdot\,\| _{H^{-1}}$ such that if $u \in H^{0}$ then
$\| u\|_{ H^{-1}}$ is the supremum of $(u,w)_{H^0}$ over all $w
\in H^1$ such that $\|w\|_{H^1} \le 1 $. $H^{-1}$ is a Hilbert
space.
\par We will write $(u,w)_{H^0}$ for $u\in H^{-1}$
and $w\in H^1$, meaning the obvious extension of the bilinear form
from $u\in H^{0}$ and $w\in H^1$.
\par
We denote by $\oo\ell _{1}$ the Lebesgue measure in $\R$, and we
denote by $ \oo{{\cal B}}_{1}$ the $\sigma$-algebra of Lebesgue sets
in $\R^1$.\
\par
 For $k=-1,0,1,2$, we  introduce the spaces
 \baaa
 X^{k}\defi L_{2}\bigl([ 0,T ], \oo{{\cal B}}_{1},\oo\ell_{1};  H^{k}\bigr),
 \quad \CC^{k}\defi C\left([0,T]; H^k\right).
  \eaaa
\par
We introduce the spaces $$ Y^{k}\defi X^{k}\!\cap \CC^{k-1}, \quad
k\ge 0, $$ with the norm $ \| u\| _{Y^k}
\defi \| u\| _{{X}^k} +\| u\| _{\CC^{k-1}}. $
\par
 We use the notations $$ \nabla u\defi \Bigl(\frac{\p u}{\p
x_1},\frac{\p u}{\p x_2},\ldots,\frac{\p u}{\p
x_n}\Bigr)^\top,\quad \nabla\cdot U=\sum_{i=1}^n\frac{\p U_i}{\p
x_i} $$ for functions $u:\R^n\to\R$ and
$U=(U_1,\ldots,U_n)^{\top}:\R^n\to \R^n$. In addition, we  use the
notation $$ (U,V)_{H^0}=\sum_{i=1}^n(U_i,V_i)_{H^0},\quad
\|U\|_{H^0}= (U,U)_{H^0}^{1/2} $$ for functions $U,V: D\to\R^n$,
where $U=(U_1,\ldots,U_n)$ and $V=(V_1,\ldots,V_n)$.
\subsubsection{The boundary value problem}  We consider the problem \be \label{4.1} \ba
\frac{\p u}{\p t} =\A u+ h,\quad \quad t\in( 0,T),\\
u(x,0)=0,\quad u(x,t)|_{x\in \p D}=0 . \ea
 \ee
 Here
 $u=u(x,t)$,
 $(x,t)\in Q$,
$h\in  L_2(Q)$, and
  \be\label{A}\A y\defi  \sum_{i=1}^n\frac{\p }{\p x_i}
\sum_{j=1}^n \Bigl(b_{ij}(x,t)\frac{\p y}{\p x_j}(x)\Bigr)
+\sum_{i=1}^n f_{i}(x,t)\frac{\p y}{\p x_i }(x)
+\,\lambda(x,t)y(x), \ee where $b(x,t):
\R^n\times[0,T]\to\R^{n\times n}$, $f(x,t):
\R^n\times[0,T]\to\R^n$, and $\lambda(x,t): \R^n\times[0,T]\to\R$,
are bounded measurable functions, and where $b_{ij}, f_i, x_i$ are
the components of $b$,$f$, and $x$. The matrix $b=b^\top$ is
symmetric.
\par
 To proceed further, we assume that Conditions
\ref{cond3.1.A}-\ref{cond3.1.B} remain in force throughout this
paper.
 \begin{condition} \label{cond3.1.A}   There exists a constant $\d>0$ such that
\be
 \label{Main1} \xi^\top  b
(x,t)\,\xi \ge \d|\xi|^2 \quad\forall\, \xi\in \R^n,\ (x,t)\in Q.
\ee
\end{condition}
\begin{condition}\label{cond3.1.B}
The functions $b(x,t):\R^n \times \R\to \R^{n \times n}$,
$f(x,t):\R^n \times \R\to \R^n$, $\lambda (x,t):\R^n \times \R\to
\R$, are differentiable in $x$, the function $b(x,t)$ is
differentiable in $t$, and the corresponding derivatives are
bounded.
\end{condition}
\par
We introduce the sets of parameters   \baaa
 &&\mu
\defi (T,\, n,\,\, D,\,\,
b(\cdot), f(\cdot),\lambda(\cdot)),
 \\
 &&{\cal P}={\cal P}(\mu)
\defi \biggl(T,\, n,\,\, D,\,\,   \delta,\,\,\,\,
\esssup_{(x,t)\in Q}\Bigl(| b(x,t)|+ |f(x,t)|+ |
\lambda(x,t)|+\Bigl|\frac{\p b}{\p x}(x,t)\Bigr|+ \Bigl|\frac{\p
f}{\p x}(x,t)\Bigr|\\ &&\hspace{10cm}+\Bigl|\frac{\p \lambda}{\p
x}(x,t)\Bigr|+\Bigl|\frac{\p b}{\p t}(x,t)\Bigr|\Bigr) \biggr).
\eaaa
 \section{Special estimates for the solution}
By the classical solvability results, there exists a unique
solution $u\in Y^1$ of problem (\ref{4.1}) for any $h\in L_2(Q)$
(see, e.g., Ladyzhenskaia (1985), Chapter III, \S \S 4-5).
Moreover,   the second energy inequality gives that, for any $K\in
\R$ and $M\ge 0$, there exists a constant $\ww C(K,M,{\cal P})>0$
such that \baa &&e^{-2Kt} (\nabla u(\cdot,t), b(\cdot,t)\nabla
u(\cdot,t))_{H^0}\nonumber\\&&
+M\biggl[e^{-2Kt}\|u(\cdot,t)\|^2_{H^0} +\int_0^t e^{-2Ks}(\nabla
u(\cdot,s), b(\cdot,s)\nabla u(\cdot,s))_{H^0}ds\biggr]
\nonumber\\ &&\hphantom{xxxxxxxxx}\le\ww C(K,M,{\cal P})\int_0^t
e^{-2Ks}\|h(\cdot,s)\|_{H^0}^2ds\qquad \forall h\in L_2(Q),\ t\in
(0,T],\hphantom{xx} \label{esti}\eaa for any solution $u$  of
problem (\ref{4.1}). This estimate follows immediately from
inequality (4.56) from Ladyzhenskaia (1985), Chapter III, and from
obvious estimates \baaa e^{-2Kt}(\nabla w, b(\cdot,t)\nabla
w)_{H^0} \le c_1\|w\|^2_{ H^1},\quad \|w\|^2_{H^0}\le
c_2e^{-2Kt}\|w\|^2_{ H^0}\quad \forall t\in[0,T], \eaaa for any
$w\in H^1$ and $w\in H^0$ respectively, where $c_1\defi
\esssup_{x,t}|b(x,t)|$ and $c_2\defi e^{2KT}$.
\par Let
$C(K,M,{\cal P})\defi\inf \ww C(K,M,{\cal P})$,  where the infimum
is taken over  all  $\ww C(K,M,{\cal P})$ such that (\ref{esti})
holds.
\begin{theorem}\label{ThM}  $$\sup_{\mu,\,M\ge 0} \inf_{K\ge 0} C(K,M,{\cal P}(\mu))\le \frac{1}{2}.
$$
\end{theorem}
\begin{corollary}\label{corr1}
For any $\mu$ and $\e>0$, there exists $ K= K(\e,{\cal P(\mu)})\ge
0$ such that \baaa \sup_{s\in[0,t]}e^{-2Ks} (\nabla u(\cdot,s),
b(\cdot,s)\nabla u(\cdot,s))_{H^0} \le \Bigl(\frac{1}{2}+\e\Bigr)
\int_0^t e^{-2Ks}\|h(\cdot,s)\|_{H^0}^2ds\\\qquad \forall h\in
L_2(Q),\ t\in(0,T], \label{sup}\eaaa where $u$ is the solution
 of problem  (\ref{4.1}) for the corresponding $h$.
\end{corollary}
\begin{remark}  For the case of a non-divergent
operator $\A u= b_{ij}u_{x_ix_j}+f_iu_{x_i}+\lambda u$, the
differentiability of the coefficients in Condition \ref{cond3.1.B}
can be weakened using the approach from Dokuchaev (2005).
\end{remark}
\section{Applications: asymptotic estimate at initial time}
 Let $$X^0_c\defi\Bigl\{h\in X^0:\ \lim_{t\to
0+}\frac{1}{t}\int_0^t\|h(\cdot,s)\|^2_{H^0}ds
=\|h(\cdot,0)\|^2_{H^0}\Bigr\}.$$ The condition that $h\in X^0_c$
is not restrictive for  $h\in X^0$; for instance, it holds if
$t=0$ is a Lebesgue point for $\|h(\cdot,t)\|_{H^0}^2$.
\begin{theorem}\label{ThAs} For any admissible set of parameters $\mu$,
 \baa\ulim_{t\to 0+}\sup_{h\in X_c^0}\frac{1}{t}\frac{
(\nabla u(\cdot,t), b(\cdot,t)\nabla u(\cdot,t))_{H^0}
}{\|h(\cdot,0)\|^2_{H^0}}\le \frac{1}{2}\,, \label{as1}\eaa where
$u$ is the solution
 of problem  (\ref{4.1}) for the corresponding $h$.
\end{theorem}
\section{On the sharpness of the estimates}
\begin{theorem}\label{lemma1} There exists a set of parameters $(n,D,
b(\cdot), f(\cdot),\lambda(\cdot))$  such that, for any $T>0$,
$M\ge 0$, \be \inf_{K\ge 0} C(K,M,{\cal P}(\mu))= \frac{1}{2}.
\label{prot1}\ee for $\mu=(T,n,D,b(\cdot),
f(\cdot),\lambda(\cdot))$.
\end{theorem}
\begin{theorem}\label{ThSharp} There exists  a set of parameters $(n,D,
b(\cdot), f(\cdot),\lambda(\cdot))$   such that \baa \ulim_{t\to
0+}\sup_{h\in X_c^0}\frac{1}{t}\frac{ (\nabla u(\cdot,t),
b(\cdot,t)\nabla u(\cdot,t))_{H^0} }{\|h(\cdot,0)\|^2_{H^0}}=
\frac{1}{2}\,, \label{as2}\eaa where $u$ is the solution
 of problem  (\ref{4.1}) for the corresponding $h$.
 \end{theorem}
  \section{Applications: parabolic equations with time delay}
 Theorem  \ref{ThM} can be also applied to analysis of parabolic
 equations with time delay. These equations have many applications, and they were intensively
 studied,
 including equations with delay
operators of a general form
 defined on the
 path of past values
 (see, e.g.,
B\'atkai and Piazzera (2001), Pao (1997), Poorkarimi and Wiener
(1999),
 Stein {\it et al} (2005), and references here). We use Theorem \ref{ThM} to obtain
 sufficient conditions of
  solvability in $Y^2$
  for the
 special case  when the first derivatives of solutions are
 affected by time delay represented by a general measurable function.
 As far as we know, this case was not
covered in
 the existing literature.
\par
Let $\tau(\cdot):[0,T]\to\R$ be a given measurable function such
that $\tau(t)\in [0,t]$.  Non-monotonic, piecewise constant, or
nowhere continuous functions $\tau(\cdot)$ are not excluded.
\par
 For functions $u:Q\to\R$, we introduce the
following operator \be\label{B1} Bu\defi\beta(x,\tau(t))^\top
\nabla u(x,\tau(t)) +\oo \beta(x,\tau(t))\,u(x,\tau(t)). \ee We
assume that the functions  $\b: \R^n\times[0,T]\to\R^n$ and
$\oo\b:\R^n\times[0,T]\to\R$  are bounded and measurable. \par Let
us consider the following boundary value problem in $Q$: \be
\label{D} \ba \frac{\p u}{\p t} =\A u+ Bu+h,\quad \quad t\in(
0,T),\\ u(x,0)=0,\quad u(x,t)|_{x\in \p D}=0 . \ea
 \ee
 Here $\A$ is such as defined above.
\par
\begin{theorem}\label{ThDelay}  Assume that there exists a constant $\d_1>0$ such that
\be
 \label{MainD} T|\xi^\top\b(x,t)|^2 \le (2-\d_1) \xi^\top  b
(x,t)\,\xi  \quad\forall\, \xi\in \R^n,\ (x,t)\in D\times [0,T].
\ee Then there exists a unique solution $u\in Y^2$ of problem
(\ref{D}) for any $h\in L_2(Q)$, and $\|u\|_{Y^2}\le
c\|h\|_{L_2(Q)}$, where $c$ is a constant that depends only on
${\cal P}$, $\d_1$, and $\sup_{x,t}|\oo\b(x,t)|$.
\end{theorem}
\begin{theorem}\label{ThDelay2}
 Assume that there exists $\t\in[0,T)$ such that $\tau(t)=0$ for
$t<\t$. Assume that the function $\tau(\cdot):[\t,T]\to \R$ is
non-decreasing and absolutely continuous, and
 $\d_*\defi \esssup_{t\in[\t,T]}\left|\frac{d\tau}{dt}(t)\right|^{-1}<+\infty$.
Then there exists a unique solution $u\in Y^2$ of problem
(\ref{D})  for any $h\in L_2(Q)$, and $\|u\|_{Y^2}\le
c_*\|h\|_{L_2(Q)}$, where $c_*$ is a constant that depends only on
${\cal P}$, $\d_*$,  $\sup_{x,t}|\b(x,t)|$, and
$\sup_{x,t}|\oo\b(x,t)|$.
\end{theorem}
\section{Proofs}
\begin{lemma}\label{lemma<<1} For any
admissible $\mu$ and any $\e>0$, $M>0$, there exists $\ww K=\ww
K(\e,M,{\cal P(\mu)})\ge 0$ such that \baa (\nabla u(\cdot,t),
b(\cdot,t)\nabla u(\cdot,t))_{H^0}+
M\biggl(\|u(\cdot,t)\|^2_{H^0}+\int_0^t (\nabla u(\cdot,s),
b(\cdot,s)\nabla u(\cdot,s))_{H^0}ds\biggr)\nonumber\\
\le\Bigl(\frac{1}{2}+\e\Bigr)\int_0^t\|h(\cdot,s)\|_{H^0}^2ds \eaa
for all $K\ge \ww K(\e,M,{\cal P})$, $t\in(0,T]$, and $h\in
L_2(Q)$, where $u\in Y^2$  is the solution of the boundary value
problem
 \be\ba \frac{\p u}{\p
t}=\A u-Ku+h,\quad t\in(0,T),\\ u(x,0) =0,\quad u(x,t)|_{x\in \p
D}=0. \label{oog3}\ea\ee
\end{lemma}
\par
Uniqueness and existence of solution  $u\in Y^2$ of problem
(\ref{oog3}) follows from the classical results (see, e.g.,
Ladyzhenskaia (1985), Chapter III).
\par {\it Proof of
Lemma \ref{lemma<<1}.} Clearly, $\A=\A_s+\A_r$, where \baaa \A_s
u=\sum_{i=1}^n\frac{\p }{\p x_i} \sum_{j=1}^n \Bigl(b_{ij}\frac{\p
u}{\p x_j}\Bigr)=\nabla\cdot (b\nabla u),\quad \quad\A_r
u=\sum_{i=1}^n f_i\frac{\p u}{\p x_i}+ \lambda u.
 \label{vv}\eaaa
Assume that the function $h(\cdot,t):D\to\R$ is differentiable and
has a compact support inside $D$ for all $t$.
 We have that \baa &&
(\nabla u(\cdot,t), b(\cdot,t)\nabla u(\cdot,t))_{H^0} -(\nabla
u(\cdot,0), b(\cdot,0)\nabla u(\cdot,0))_{H^0}\nonumber\\&& =
2\int_0^t \left(\nabla u,b\nabla\frac{\p u}{\p s}\right)_{H^0} ds
+ \int_0^t \left(\nabla u,\frac{\p b}{\p s}\nabla u \right)_{H^0}
ds \nonumber
\\&&=2\int_0^t \left(\nabla u,b\nabla\bigl(\A
u-Ku+h\bigr)\right)_{H^0} ds  + \int_0^t \left(\nabla u,\frac{\p
b}{\p s}\nabla u \right)_{H^0} ds \nonumber\\&& =2\int_0^t
\left(\nabla u,b\nabla\bigl(\nabla\cdot (b\nabla
u)\bigr)\right)_{H^0} ds +2\int_0^t \left(\nabla u,b\nabla\A_r
u\right)_{H^0} ds\nonumber\\&&\hspace{0.3cm}-2K\int_0^t (\nabla u,
b\nabla u)_{H^0} ds+2\int_0^t \left(\nabla u,b\nabla
h\right)_{H^0} ds +\int_0^t \left(\nabla u,\frac{\p b}{\p s}
\nabla u \right)_{H^0} ds.\hphantom{xx} \label{R0} \eaa
\par Let an arbitrary $\e_0>0$ be
given.  We have that in under the integrals in (\ref{R0}),
 \baa \left(\nabla u, \frac{\p b}{\p t} \nabla
u \right)_{H^0}\le  \esssup_{x,t} \left|\frac{\p b}{\p
t}(x,t)\right| \|\nabla u\|^2_{H^0} \le c_{\d}' (\nabla
u(\cdot,t), b(\cdot,t)\nabla u(\cdot,t))_{H^0}, \label{dv}\eaa
where $c_{\d}'=c_{\d}'({\cal P})$ is a constant that depends on
${\cal P}$ only. Further,
 \baa 2\left(\nabla u,b\nabla h\right)_{H^0} =-2\left(\nabla\cdot (b\nabla u), h\right)_{H^0}\le
\frac{2}{1+2\e_0}(\nabla u, b \nabla
u)^2_{H^0}+\left(\frac{1}{2}+\e_0\right)\left\|h\right\|^2_{H^0},
\label{R1}\eaa
\par
Let a bounded measurable function $v(x,t):Q\to\R$ be such that
$b=v^\top v$. We have that
 \baaa 2\left(\nabla u,b\nabla \A_r u\right)_{H^0} = 2\left(v\nabla u,v\nabla \A_r u\right)_{H^0} \le
\e_1^{-1}\left\|v\nabla
u\right\|^2_{H^0}+\e_1\left\|v\nabla\A_ru\right\|^2_{H^0}\quad\forall\e_1>0.
\label{R3'}\eaaa It follows that
 \baaa
2\left(\nabla u,b\nabla \A_r u\right)_{H^0} \le \e_1^{-1}(\nabla
u,b \nabla u)_{H^0}+c_1\e_1\left\|u\right\|^2_{H^2}, \eaaa where
$c_1=c_1({\cal P})$ is a constant that depends on ${\cal P}$ only.
By the second fundamental inequality, there exists a constant
$c_{*}=c_{*}({\cal P})>0$ such that
\be\int_0^t\left\|u(\cdot,s)\right\|^2_{H^2}ds\le
c_{*}\int_0^t\left\|h(\cdot,s)\right\|^2_{H^0}ds. \label{e*} \ee
(See, e.g., estimate (4.56) from Ladyzhenskaia (1985), Chapter
III). By Lemma 5.3 from Dokuchaev (2005), one can choose the same
constant $c_*$  for all $t\in [0,T]$, $K>0$. Hence \baa
2\int_0^t\left(\nabla u,b\nabla \A_r u\right)_{H^0}ds \le
\e_1^{-1}\int_0^t(\nabla u,b \nabla u)_{H^0}
ds+\e_0\int_0^t\left\|h\right\|^2_{H^0}ds, \label{R3} \eaa
 where $\e_1>0$ is such that $c_1c_{*}\e_1=\e_0$.
\par
By Lemma 5.2 from Dokuchaev (2005), p. 357, it follows that there
exists $K_1=K_1(\e_0,M,{\cal P})>0$ such that if $K>K_1$, then \baa
M\sup_{s\in[0,t]}\left\|u(\cdot,s)\right\|^2_{H^0}\le \e_0\int_0^t
\|h(\cdot,s)\|_{H^0}^2ds. \label{R4}\eaa \par
\par
 Assume  now  that $D=\R^n$. In
that case,  we have immediately that \baa 2\left(\nabla
u,b\nabla\bigl(\nabla\cdot (b\nabla
u)\bigr)\right)_{H^0}=2\left(b\nabla
u,\nabla\bigl(\nabla\cdot(b\nabla u)\bigr)\right)_{H^0}
=-2\left(\nabla\cdot (b\nabla u ),\nabla\cdot(b\nabla
u)\right)_{H^0}\nonumber\\=-2\left\|\nabla\cdot(b\nabla u
)\right\|^2_{H^0}.\label{R2}\eaa The third equality here was
obtained using integration by parts.
\par
 By (\ref{R0})-(\ref{R4}), it
follows that if $K>K_1$ and $2K>\e_1^{-1}+c_{\d}'+M$, then \baa &&
(\nabla u(\cdot,t), b(\cdot,t)\nabla
u(\cdot,t))_{H^0}+M\biggl(\|u(\cdot,t)\|^2_{H^0}+\int_0^t (\nabla
u(\cdot,s), b(\cdot,s)\nabla u(\cdot,s))_{H^0} ds\biggr)
\nonumber\\&&\le \Bigl(\frac{2}{1+2\e_0}-2\Bigr)\int_0^t
\left\|\nabla\cdot(b\nabla u )\right\|^2_{H^0} ds
+\left(\e_1^{-1}+c_{\d}'+M-2K\right)\int_0^t \left\|v\nabla
u\right\|^2_{H^0} ds \nonumber
\\&&
+\Bigl(\frac{1}{2}+3\e_0\Bigr)\int_{0}^t\left\|h(\cdot,s)\right\|^2_{H^0}ds
\le
\Bigl(\frac{1}{2}+3\e_0\Bigr)\int_0^t\|h(\cdot,s)\|_{H^0}^2ds.\hphantom{}\label{R4'}
\eaa  Then the proof of Lemma \ref{lemma<<1} follows for the case
when $D=\R^n$, since (\ref{R4'}) holds for all $h$ from a set that
is  dense in $L_2(Q)$.
\par
To complete the proof of Lemma \ref{lemma<<1}, we need to cover
the case  when $D\neq \R^n$. From now and up to the end of the
proof of this lemma, we assume that $D\neq \R^n$. In that case,
(\ref{R2}) does not hold, since the integration by parts used for
the third inequality in (\ref{R2}) is not applicable anymore. To
replace (\ref{R2}), we are going to show that there exists a
constant
 $C=C({\cal P})>0$ such that, for an arbitrarily $\e_2>0$,
\baa \left(\nabla u,b\nabla\bigl(\nabla\cdot (b\nabla
u)\bigr)\right)_{H^0}\le -\left\|\nabla\cdot(b\nabla u
)\right\|^2_{H^0}  +\e_2 \|u\|^2_{H^2}+C\e_2^{-1}\|u\|^2_{H^1}.
\label{estD}\eaa
\par
 Assume that (\ref{estD}) holds.  Since $D$ is bounded, we have that
 $\|u\|^2_{H^1}\le c_{\d}(\nabla u,b\nabla u)^2_{H^0}$ for some constant
 $c_{\d}= c_{\d}({\cal P})$.  In addition,  (\ref{R1}) is still valid, since we assumed that
 $h(\cdot,t)$ have support inside $D$. Similarly
 to (\ref{R4'}), we obtain that if  $2K>\e_1^{-1}+c_{\d}'+Cc_{\d}\e_2^{-1}+M$ and
 $K>K_1$, then
\baaa &&(\nabla u(\cdot,t), b(\cdot,t)\nabla u(\cdot,t))_{H^0} +
M\biggl[\|u(\cdot,t)\|^2_{H^0}+\int_0^t (\nabla u(\cdot,s),
b(\cdot,s)\nabla u(\cdot,s))_{H^0} ds\biggr] \\&&\le
\Bigl[\frac{2}{1+2\e_0}-2\Bigr]\int_0^t \left\|\nabla\cdot
(b\nabla u)\right\|^2_{H^0} ds+[\e_1^{-1}+c_{\d}'+M-2K]\int_0^t
(\nabla u(\cdot,s), b(\cdot,s)\nabla u(\cdot,s))_{H^0} ds
\\&&+
\Bigl(\frac{1}{2}+3\e_0\Bigr)\int_{0}^t\left\|h(\cdot,s)\right\|^2_{H^0}ds
+2\e_2
\int_{0}^t\|u(\cdot,s)\|^2_{H^2}ds+2C\e_2^{-1}\int_{0}^t\|u(\cdot,s)\|^2_{H^1}ds
\\&&\hspace{5cm}\le
\Bigl(\frac{1}{2}+3\e_0+2c_{*}\e_2\Bigr)\int_0^t\|h(\cdot,s)\|_{H^0}^2ds,\hphantom{}\label{R4''}
\eaaa where $c_{*}$ is the constant from (\ref{e*}). Then the
proof of Lemma \ref{lemma<<1} follows provided that  (\ref{estD})
holds.
 Therefore, it suffices to prove (\ref{estD}) for $D\neq \R^n$.
\par
 Let us prove (\ref{estD}). This estimate can rewritten as
  \baaa \left(w,\nabla\bigl(\nabla\cdot w \bigr)\right)_{H^0} \le
-\left((\nabla\cdot w),\bigl(\nabla\cdot w \bigr)\right)_{H^0}
+\e_2 \|u\|^2_{H^2}+C\e_2^{-1}\|u\|^2_{H^0}, \eaaa
 where $w=b\nabla u$.
 We have
\baa \left(w,\nabla\bigl(\nabla\cdot w \bigr)\right)_{H^0} =
-\left((\nabla,w),\bigl(\nabla\cdot w \bigr)\right)_{H^0}
+\sum_{i=1}^n\int_{\p D}\w J_i dz.\label{estD1} \eaa Here  $\w
J_{i} = z(\nabla,z)\cos({\bf n},e_i)$, where ${\bf n}={\bf n}(s)$
is the outward pointing normal to the surface $\p D$  at the point
$s \in \p D$, and $e_k=(0,\ldots,0,1,0,\ldots,0)$ is the $k$th
basis vector in the Euclidean space $\R^n$.  We have that \baaa\w
J_{i} =\sum_{j,k,m=1}^n\a_{ijkm}
J_{ijkm}+\sum_{j,k,m=1}^n\a'_{ijkm} J'_{ijkm}, \eaaa where \baaa
J_{ijkm} = \frac{\p u}{\p x_j}\frac{\p^2 u}{\p x_k x_m}\cos({\bf
n},e_i),\quad J'_{ijk} = \frac{\p u}{\p x_j}\frac{\p u}{\p
x_k}\cos({\bf n},e_i),\eaaa and where $\a_{ijkm}$, $\a_{ijk}'$ are
some bounded functions.
\par
Let us estimate $\int_{\p D}\w J_{i}dz$.  We mostly follow the
approach  from Ladyzhenskaya and Ural'tseva (1968), Section 3.8.
Let $x^0=\{x_i^0\}^n_{i=1} \in \p D$ be an arbitrary point. In its
neighborhood, we introduce local Cartesian coordinates
$y_m=\sum_{k=1}^{n}c_{mk}(x_k-x^0_k)$ such that the axis  $y_n$ is
directed along the outward normal ${\bf n}={\bf n}(x_0)$ and
$\{c_{mk}\}$ is an orthogonal matrix.
\par
Let $y_n=\psi (y_1,\ldots,y_{n-1})$ be an equation determining the
surface  $\p D$ in a neighborhood of the origin. By the properties
of the surface $\p D$, the first order and second order
derivatives of the function  $\psi $ are bounded. Since
$\{c_{mk}\}$ is an orthogonal matrix, we have
$x_k-x^0_k=\sum_{m=1}^{n}c_{km}y_m$.  Therefore, $\cos({\bf
n},e_m)=c_{nm}$, $m=1,\ldots,n$. Then \baaa  J_{ijkm}
=\sum_{l=1}^nc_{jl}\frac{\p u}{\p y_l}
\sum_{p,q=1}^{n}c_{pk}c_{qm} \frac{\p^2 u}{\p y_p\p y_q}
c_{ni},\quad J'_{ijk} =\sum_{l=1}^nc_{jl}\frac{\p u}{\p y_l}
\sum_{p=1}^{n}c_{pk} \frac{\p u}{\p y_p} c_{ni}. \eaaa The
boundary condition $u(x,t)|_{x\in \p D}=0$ can be rewritten as  $$
     u(y_1,\ldots,y_{n-1},\psi (y_1,\ldots,y^{n-1}),t)=0
$$ identically  with respect to $y_1,\ldots,y_{n-1}$ near the
point $y_1=\ldots=y_{n-1}=0$. Let us differentiate this identity
with respect to  $y_p$ and $y_q$, $p,q=1,\ldots,n-1$, and take
into account that $$ \frac{\p \psi}{\p y_p}=0, \quad
p=1,\ldots,n-1. $$ at $x_0$. Then $$ \frac{\p u}{\p y_p}=0, \quad
\frac{\p^2 u}{\p y_p\p y_q}  =- \frac{\p u}{\p y_n} \frac{\p^2
\psi}{\p y_p\p y_q} = -\frac{\p u}{\p {\bf n}} \frac{\p^2 \psi}{\p
y_p\p y_q}, \quad p,q=1,\ldots,n-1. $$ Hence \baa \int_{\p D}\w
J_{i}(z,s)dz =\int_{\p D}\sum_{j,k,m=1}^n\a_{ijkm}
J_{ijkm}(z,s)dz+\sum_{j,k,m=1}^n\int_{\p D}\a'_{ijkm} J'_{ijkm}
\a_{ijkm}(z,s)dz \nonumber\\ \hphantom{}\le  \w c_1\int_{\p D}
\biggl| \frac{\p u}{\p {\bf n}}\biggr|^2 dz\le \e_2\sum_{i,j=1}^n
\int_{ D} \biggl|\frac{\p^2 u}{\p x_i\p x_j}(x,s)\biggr|^2dx +\w
c_2(1+\e^{-1}_2) \|u(\cdot,s)\|^2_{H^1} \quad \forall \e_2>0
\label{3.3.12} \eaa for some constants $\w c_i=\w c_i({\cal P})$.
The last estimate follows from the estimate (2.38) from
Ladyzhenskaya and Ural'tseva (1968), Chapter II. By (\ref{estD1})
and (\ref{3.3.12}), it follows (\ref{estD}). This completes the
proof of Lemma \ref{lemma<<1}. $\Box$
\par
{\it Proof of Theorem \ref{ThM}}. Clearly,
$u(x,t)=e^{Kt}u_K(x,t)$, where $u$ is the solution
 of problem  (\ref{4.1})  and $u_K$ is the solution of (\ref{oog3}) for
the  nonhomogeneous  term $e^{-Kt}h(x,t)$. Therefore, Theorem
\ref{ThM} follows immediately from Lemma \ref{lemma<<1}. $\Box$
 \par  Corollary \ref{corr1}
 follows immediately from Theorem \ref{ThM}.
\par
{\it Proof of Theorem \ref{ThAs}.} Let $\e>0$ be given. By
Corollary \ref{corr1}, there exists $K(\e)=K(\e,{\cal P}(\mu))$
such that \baa e^{-2K(\e)t}(\nabla u(\cdot,t), b(\cdot,t)\nabla
u(\cdot,t))_{H^0}\le
\left(\frac{1}{2}+\e\right)\int_0^te^{-2K(\e)s}\|h(\cdot,s)\|^2_{H^0}ds\nonumber\\
\qquad \forall t\in(0,T),\ h\in X^0. \label{aas}\eaa
\par
Let $p(h,t)\defi\frac{1}{t}\int_0^t\|h(\cdot,s)\|^2_{H^0}ds$ and
$q(u,t)\defi (\nabla u(\cdot,t), b(\cdot,t)\nabla
u(\cdot,t))_{H^0}$.  It follows that
 \baaa
\sup_{h\in
X^0}\left(\frac{q(u,t)}{tp(h,t)}-\frac{1-e^{-2K(\e)t}}{tp(h,t)}q(u,t)
\right)\le \left(\frac{1}{2}+\e\right)
 \qquad \forall t\in(0,T). \eaaa
 Hence
 \baaa
\sup_{h\in X^0}\frac{1}{tp(h,t)}q(u,t)\le
\left(\frac{1}{2}+\e\right)+\sup_{h\in
X^0}\frac{1-e^{-2K(\e)t}}{tp(h,t)}q(u,t).
 \qquad \forall t\in(0,T). \eaaa
By (\ref{aas}),
 \baaa q(u,t)\le
e^{2K(\e)t}\left(\frac{1}{2}+\e\right)tp(h,t) \qquad \forall
t\in(0,T),\ h\in X^0. \eaaa Hence \baaa\sup_{h\in
X^0}\frac{1-e^{-2K(\e)t}}{tp(h,t)}q(u,t)\to 0 \quad\hbox{as}\quad
t\to 0+\quad \forall  \e>0. \eaaa
\par
If $h\in X^0_c$, then $p(h,t)\to \|h(\cdot,0)\|^2_{H^0}$ as $t\to
0+$. It follows that \baaa\ulim_{t\to 0+}\sup_{h\in
X^c}\frac{q(u,t)}{t\|h(\cdot,0)\|^2_{H^0}}\le
\left(\frac{1}{2}+\e\right)
 \eaaa
for any $\e>0$. Then (\ref{as1}) follows.
 This completes the proof of Theorem \ref{ThAs}.
 $\Box$
\par
{\it Proof of Theorem \ref{lemma1}.}
 Repeat that
$u(x,t)=e^{Kt}u_K(x,t)$, where $u$ is the solution
 of problem  (\ref{4.1})  and $u_K$ is the solution of (\ref{oog3}) for
$h_K(x,t)=e^{-Kt}h(x,t)$. Therefore,  it suffices to find $n$,
$D$, $b,f,\lambda$  such that\baa &&\forall T>0,c>0, K>0\quad
\exists h\in L_2(Q):\quad  \nonumber\\   &&(\nabla u(\cdot,T),
b(\cdot,T)\nabla u(\cdot,T))_{H^0}
\ge\left(\frac{1}{2}-c\right)\int_0^T\|h(\cdot,t)\|_{H^0}^2dt,
\label{prot10}\eaa  where $u$ is the solution
 of problem  (\ref{oog3}).
\par
 Let us show that (\ref{prot10}) holds for \be
 n=1,\quad D=(-\pi,\pi),\quad b(x,t)\equiv 1,
 \quad f(x,t)\equiv 0,\quad \lambda(x,t)\equiv 0.\label{mu+}\ee
\par
 In this case,
(\ref{oog3}) has the form \baa u'_t=u''_{xx}-Ku+h,\quad
u(x,0)\equiv 0,\quad u(x,t)|_{x\in\p D}=0, \label{eu}\eaa
\par
Let
 \baa
 \g=m^2+K,\quad
h_m(x,t)\defi\sin(m x) e^{\g t}, \quad \label{wxi}\eaa where
$m=1,2,3,\ldots$.
 It can be verified immediately that  the
solution of the boundary value problem is \baaa
u(x,t)=\sin(mx)\int_0^te^{-\g(t-s)+\g s}ds=\sin(m x)e^{-\g
t}\int_0^te^{2\g s}ds =\sin(m x)e^{-\g t}\frac{e^{2\g t}-1}{2\g }.
\eaaa Hence \baaa\|u'_x(\cdot,T)\|^2_{H^0}=m^2 \|\cos(m
x)\|^2_{H^0}e^{-2\g T}\left(\frac{e^{2\g }-1}{2\g}\right)^2
=m^2\pi e^{-2\g T}\frac{(e^{2\g T}-1)^2}{4\g^2}, \eaaa and \baaa
\int_0^T\|h(\cdot,t)\|^2_{H^0}dt=\|\sin(m
x)\|^2_{H^0}\int_0^Te^{2\g t}dt= \pi\frac{e^{2\g T}-1}{2\g}. \eaaa
It follows that \baa
\|u'_x(\cdot,T)\|^2_{H^0}\left(\int_0^T\|h(\cdot,t)\|^2_{H^0}dt\right)^{-1}
=\frac{m^2}{2\g}e^{-2\g T}(e^{2\g T} -1)=\frac{m^2}{2\g}(1-e^{-2\g
T})\to \frac{1}{2} \label{1/4}\eaa as $\g\to +\infty$. In
particular, it holds if $K$ is fixed and $m\to +\infty$. It
follows that (\ref{prot1}) holds.
 This completes the proof of Theorem \ref{lemma1}. $\Box$
\par
{\it Proof of Theorem \ref{ThSharp}.} Let the parameters of the
equation be defined by (\ref{mu+}). Consider a sequence $\{T_i\}$
such that $T_i\to 0+$ as $i\to +\infty$.  Let $h=h_{m}$ be defined
by (\ref{wxi}) for an increasing sequence of integers $m=m_i$ such
that $m_i>T_i^{-1}$. In that case,
(\ref{1/4}) holds since $\g T\to +\infty$.
 Then (\ref{as2})
holds.  $\Box$
\par {\it Proof of Theorem
\ref{ThDelay}.} For $K>0$, introduce operators \be\label{BK}
B_Ku\defi e^{K(\tau(t)-t)}\Bigl(\beta(x,\tau(t))^\top\nabla u
(x,\tau(t))+\oo \beta(x,\tau(t))\,u(x,\tau(t))\Bigr). \ee
 Note that $u\in Y^2$ is the solution of the problem (\ref{D})
 if and only if
 $u_K(x,t)=e^{-Kt}u(x,t)$ is the solution of the problem  \be \label{D2}\ba \frac{\p
u_K}{\p t} =\A u_K-Ku_K +B_Ku_K+h_K,\quad \quad t\in( 0,T),\\
u_K(x,0)=0,\quad u_K(x,t)|_{x\in \p D}=0,\ea \ee where
$h_K(x,t)=e^{-Kt}h(x,t)$.
 In addition,
 $$ \|u\|_{Y^2}\le
e^{KT}\|u_K\|_{Y_2},\quad \|h_K\|_{L_2(Q)}\le \|h\|_{L_2(Q)}. $$
Therefore, uniqueness and solvability in $Y^2$ of problem
(\ref{D}) follows from  existence of $K>0$ such that problem
(\ref{D2}) has an unique solution in $Y^2$. Let us show that this
$K$ can be found.
\par
We introduce  operators $F_K: L_2(Q)\to Y^2$ such that $u=F_Kh$ is
the solution of problem (\ref{oog3}).
\par
 Let $g\in L_2(Q)$ be such that
\baa g=h+B_Kw,\quad  \hbox{where}\quad w=F_Kg.\label{gw}\eaa It
can be rewritten as  $g=h+R_Kg$, or \baa g-R_Kg=h,\label{gw2}\eaa
where the operator $R_K:L_2(Q)\to L_2(Q)$ is defined as
$$R_K=B_KF_K.$$ In that case, $u_K\defi F_Kg\in Y^2$ is the
solution of (\ref{D2}).
\par
Let us show that there exists $K>0$ such that \baa\|R_K\|<1.
\label{<1}\eaa
\par Let $w=F_Kh$. By Theorem \ref{ThM} reformulated as  Lemma \ref{lemma<<1},
for  any $\e>0$, $M>0$, there exists $K(\e,M,{\cal P(\mu)})\ge 0$
such that \baa \sup_{t\in[0,T]} \Bigl((\nabla w(\cdot,t),
b(\cdot,t)\nabla w(\cdot,t))_{H^0}+
M\|w(\cdot,t)\|^2_{H^0}+M\int_0^t (\nabla w(\cdot,s),
b(\cdot,s)\nabla w(\cdot,s))_{H^0}ds\Bigr)\nonumber\\
\le\Bigl(\frac{1}{2}+\e\Bigr)\|h\|_{L_2(Q)}^2ds \quad \forall h\in
L_2(Q).\hphantom{xxx}\label{ineq} \eaa \par Let $\w w(x,t)\defi
w(x,\tau(t))$, $\w\b(x,t)\defi \b(x,\tau(t))$, and $\w b(x,t)\defi
b(x,\tau(t))$.
 By the
definitions, \baa \|R_Kh\|_{L_2(Q)}=\|B_Kw_K\|_{L_2(Q)}\le
\|\w\b^\top\nabla \w w\|_{L_2(Q)}+C_\b\|\w w\|_{L_2(Q)},
\label{RK}\eaa where $C_\b\defi\sup_{x,t}|\oo\b(x,t)|$. Clearly,
\baa  &&\|\w\b^\top\nabla \w w\|_{L_2(Q)}^2\le
T\sup_{t\in[0,T]}\|\w\b(\cdot,t)^\top\nabla \w
w(\cdot,t)\|^2_{H^0}\le T\sup_{t\in[0,T]}\|\b(\cdot,t)^\top\nabla
 w(\cdot,t)\|^2_{H^0},\label{RK20}\\
 &&\|\w w\|_{L_2(Q)}^2\le
T\sup_{t\in[0,T]}\|\w w(\cdot,t)\|_{H^0}^2\le
T\sup_{t\in[0,T]}\|w(\cdot,t)\|_{H^0}^2. \label{RK2}\eaa
\par
Further, let   $K>0$, $M>0$, and $\e>0$, be such that (\ref{ineq})
is satisfied and \baa \d_2<1, \quad
\d_2^2+\d_3^2+2\d_2\d_3<1,\label{delta23}\eaa where \baa\d_2\defi
\sqrt{\left(2-\d_1\right)\left(\frac{1}{2}+\e\right)},\qquad
 \quad \d_3\defi C_\b\sqrt{T
M^{-1}\Bigl(\frac{1}{2}+\e\Bigr)}. \label{delta2}\eaa
 By (\ref{ineq}), (\ref{RK20})-(\ref{delta2}),and (\ref{MainD}), it follows that \baaa\|\w\b^\top\nabla \w
w\|_{L_2(Q)}^2\le T\sup_{t\in[0,T]}\|\b(\cdot,t)^\top\nabla
w(\cdot,t)\|^2_{H^0} \le (2-\d_1)\sup_{t\in[0,T]}(\nabla  w, b
\nabla  w)_{H^0}\le \d_2^2\|h\|_{L_2(Q)}^2.  \eaaa In addition, we
have that \baaa C_\b^2\|\w w\|^2_{L_2(Q)}\le
TC_\b^2\sup_{t\in[0,T]}\| w(\cdot,t)\|_{H^0}^2 \le TC_\b^2
M^{-1}\Bigl(\frac{1}{2}+\e\Bigr)\|h\|_{L_2(Q)}^2ds\le
\d_3^2\|h\|_{L_2(Q)}^2.  \eaaa By (\ref{RK}), it follows that
\baaa
 \|R_Kh\|^2_{L_2(Q)}\le
\left(\|\w\b^\top\nabla \w w\|_{L_2(Q)}+C_\b\|\w
w\|_{L_2(Q)}\right)^2
\le
\left(\d_2^2 +2\d_2\d_3+\d_3^2\right) \|h\|^2_{L_2(Q)}.
 \eaaa
By (\ref{delta23}), it follows that (\ref{<1}) holds, where the
norm of the operator $R_K:L_2(Q)\to L_2(Q)$ is considered. It
follows that the operator
 $(I-R_K)^{-1}:L_2(Q)\to L_2(Q)$ is continuous, By (\ref{gw})-(\ref{gw2}),
$u_K=F_K(I-R_K)^{-1}h_K$ is the solution of problem (\ref{D2}) for
$h_K\in L_2(Q)$.
\par
The choice of $K$, $M$, and $\e$, depends on ${\cal P}$, $\d_1$, and
$\sup_{x,t}|\oo\b(x,t)|$ only. Hence $\d_2$ and $\d_3$ depends on
these parameters only. It follows that the norm of the operator
 $(I-R_K)^{-1}:L_2(Q)\to L_2(Q)$ can be estimated from above by a constant
 that
depends only on these parameters.  This proves the  estimate
 for the solution stated in Theorem
\ref{ThDelay}. This completes the proof of this theorem. $\Box$
\par
{\it Proof of Theorem \ref{ThDelay2}} is based again on
(\ref{ineq}). It is  similar to the proof of Theorem
\ref{ThDelay}, with a minor modification: instead of (\ref{RK20}),
we use that \baaa &&\int_0^T\|\b(\cdot,\tau(t))^\top \nabla
u(\cdot,\tau(t))\|_{H^0}^2dt=\int_{\t}^T\|\b(\cdot,\tau(t))^\top
\nabla
u(\cdot,\tau(t))\|_{H^0}^2dt\\&&=\int_{\t}^T\|\b(\cdot,\tau(t))^\top
\nabla
u(\cdot,\tau(t))\|_{H^0}^2\left(\frac{d\tau(t)}{dt}\right)^{-1}d\tau(t)
\le
\d_*\int_{\tau(\t)}^{\tau(T)}\|\b(\cdot,s)^\top \nabla
u(\cdot,s)\|_{H^0}^2ds\\&&\le c\int_{0}^{T}(\nabla u(\cdot,s),
b(\cdot,s)\nabla u(\cdot,s))_{H^0}
 ds, \eaaa where $c$
is a constant that depends on $\d$, $\d_*$, and
$\sup_{x,t}|\b(x,t)|$. $\Box$ \par It can be seen from the proofs
that the approach used for Theorems \ref{ThDelay}-\ref{ThDelay2} can
be extended on more general delay operators  represented by
integrals accumulating the past values.
\subsection*{Acknowledgment}  This work  was supported by NSERC
grant of Canada 341796-2008 to the author.
\section*{References}
$\hphantom{xx}$ B\'atkai, A., Piazzera, A.S. (2001). Semigroups
and linear partial differential equations with delay. {\it J.
Math. Anal. Appl.} {\bf 264}, 1--20.
\par
 Dokuchaev, N.G. (2005). Parabolic Ito equations
and second fundamental inequality.  {\it Stochastics} {\bf 77},
iss. 4., 349-370.
\par
 Ladyzhenskaya, O. A., and Ural'tseva, N.N.  (1968). {\it
  Linear and quasilinear elliptic equations}. New York: Academic Press.
\par
Ladyzhenskaia, O.A. (1985). {\it The Boundary Value Problems of
Mathematical Physics}. New York: Springer-Verlag.
\par
Pao, C.V. (1997).  System of Parabolic Equations with Continuous
 and Discrete Delays, {\it J. Math.
Anal. Appl.}, {\bf 205}, 157--185.
\par
Poorkarimi, C.H.,  Wiener, J. (1999).  Bounded solutions of
nonlinear parabolic equations with time delay. {\it  Electronic J.
Differential Equations},  1999, pp. 87--91.
\par
Stein, M., Vogt H., V\"oigt, J. (2005). The modulus semigroup for
linear delay equations III. {\it J. Funct. Anal.}, {\bf 220} (2),
388--400.
\end{document}